\documentclass{article}

\usepackage{amssymb}
\usepackage{amsmath}
\usepackage[dvips]{graphicx}
\begin{document}
\begin{center}
\noindent \textbf{ SPECTRAL ANALYSIS FOR ONE CLASS OF SECOND-ORDER
INDEFINITE NON-SELF-ADJOINT DIFFERENTIAL OPERATOR PENCIL }
\end{center}

\begin{center}
 \textbf{R.F. Efendiev }
\end{center}
\begin{center}
\textbf{Institute of  Applied Mathematics, Baku State University,}

\textbf{Z.Khalilov, 23, AZ1148, Baku, Azerbaijan, }
\end{center}
\begin{center}
efendievrakib@bsu.az
\end{center}

\indent{\underline{\textbf{ABSTRACT}}}

\noindent

The inverse problem for the Schrodinger operator with complex
periodic potentials and  discontinuous coefficient on the axis is
studied.  Main characteristics of the fundamental solutions are
investigated, the spectrum of the operator is studied.  We give
the formulation of the inverse problem, prove a uniqueness theorem
and provide a constructive procedure for the solution of the
inverse problem.

\textbf{\underbar{Key words}}:  Discontinuous equation, Spectral
singularities, Inverse spectral problem, Continuous spectrum.

\textbf{\underbar{MSC}}:   34A36, 34L05, 47A10, 47A70

\begin{center}
\indent{\underline{\textbf{INTRODUCTION}}}
\end{center}

    In this paper we propose a method for
     solving the inverse spectral problem of wave propagation with
     discontinuous wave speed in one-dimensional layered-inhomogeneous
      medium in the frequence domain which is described by
      Schrodinger equation
 \begin{equation} \label{GrindEQ__1_}
-y''(x)+2\lambda p(x)y(x)+q(x)y(x)=\lambda ^{2} \rho
\left(x\right)y(x)
\end{equation}
where $R$- real axis , $\lambda$- wave number ( momentum),
$\lambda^{2}$- energy, $p(x)$- describes the joint effect of
absorption and generation of energy, $q(x)$- describes the
regeneration of the force density , $\rho^{-1}(x)$-  wave speed.

\indent   An inverse problem consists in inferring some properties
of given problem from a complete or partial knowledge of its
solution. The definition is rather vague: nevertheless it points
out that the difference between "direct problem" and an "inverse
problem" is, to some extent, arbitrary and its rational is to be
found in  the historical development of that particular problem.
At this time, there is no general inverse spectral theory-even the
simplest cases require considerable ingenuity for their resolution
and none of inverse problems can be really posed
 unless a class of coefficients  is specified in advance.

\indent    The direct and inverse scattering problems for the
classical Schrodinger operators  have been studied extensively
(see [ 1-7] and references therein).

 \indent  For the Sturm-Liouville operator with different singularities (i.e. on the
half-line  having an arbitrary number of turning points, having
singularities and turning points at   the end-points of the
interval) the determination of the spectral function or
  normalizing constants has been studied by G. Freiling and V.Yurko [8-11]
 These results are mainly based on Volterra operator transformation and contour integrations.

\indent The string equation or wave equation in layered medium was
investigated by  M.G. Krein [12],A.S. Blagoveshchensky [13],
  N.I.Grinberg [14],   and others.

\indent M.I.Belishev [15] first considered (1) with
$p(x)=0,q(x)=0$ and solved the inverse problem of reconstructing
the $\rho(x)$ of a finite inhomogeneous string from the
frequencies and energies of its normalized characteristic
vibrations in the case where $\rho(x)$ can change sign (is
indefinite).

\indent  The inverse spectral problem of recovering pencils of
second-order differential operators on the half-line with turning
points was studied by V.Yurko in [12],where it was established
properties of spectral characteristics, gave  formulation of the
inverse problem and prove a uniqueness theorem for solution of the
inverse problem.

\indent The considered problem is of interest not only for its own
sake,  but also because there are other inverse problems in
physics   which can be reduced to it. For example, kinetic
equation  after separation of variables reduce Sturm-Liouville
eigenvalue problems with an indefinite function.

\indent  In neutron transport one may deal with groups of neutrons
of different speeds, anisotropic interactions and interactions
between different groups of neutrons, which in general leads to
nonselfadjoint  operator.

 \indent In present work we suggested that layered-inhomogeneous medium
  can also absorption and emit an energy. Therefore, the potentials
  $p (x),q(x)$ are complex valued.

   \indent Complex potentials have been used in the past to model
neutron absorption [16,17]  and thus the imaginary potential
represents emission or absorption. For example, in the crystal
model the complex potential is viewed as representing a "prepared"
material. The material is assumed to be capable of absorption  or
emission.

\indent Now we begin to study  the inverse problem for a pencil
$L$ of non-self-adjoint differential operators pencil  generated
by formal differential expression

\begin{equation}\label{GrindEQ__2}
 l(\frac{d}{{dx}},\lambda ) \equiv \frac{1}{{\rho \left( x
\right)}}\left\{ { - \frac{{d^2 }}{{dx^2 }} + 2\lambda p\left( x
\right) + q\left( x \right)} \right\}
\end{equation}

in the space $L_{2}(-\infty,\infty,\rho(x))$, $ \lambda$  is
complex, and the coefficients $p(x)$, $q(x)$  and $\rho(x)$ are
\begin{equation} \label{GrindEQ__3_}
p(x)=\sum _{n=1}^{\infty }p_{n} e^{inx}  ;  \sum _{n=1}^{\infty
}n\cdot \left|p_{n} \right|<\infty
\end{equation}

\begin{equation}\label{GrindEQ__4_}
q(x)=\sum _{n=1}^{\infty }q_{n} e^{inx}  ; \sum _{n=1}^{\infty
}\left|q_{n} \right|<\infty \, \, \, \, \, \,
\end{equation}

 \begin{equation}
\label{GrindEQ__5_} \rho \left(x\right)=\left\{\, \,
\begin{array}{c} {1} \\ {-\beta ^{2} } \end{array}\right. \, \, \,
\begin{array}{c} {for} \\ {for} \end{array}\, \, \, \,
\begin{array}{c} {x\geq 0,} \\ {x<0.\, \, } \end{array}
\end{equation}

\indent Indefinite differential pencils (1) with coefficients
(3-5),  produce qualitative modifications in the investigations of
the inverse problem.

The feature of our case is the fact that an operator of
transformation is absent and because the  function $\rho
\left(x\right)$ is nonsmooth, the  wave equation is  cannot be
reduced to the Schrodinger equation by  change of variables.

\indent Firstly potential \eqref{GrindEQ__4_} was considered by
M.G.Gasymov[18].  As a final remark, we mention some other works
related to the potentials
\eqref{GrindEQ__3_}-\eqref{GrindEQ__4_}[20-23].

\indent Note that at $p(x)=0$ given problem was investigated
 in the work [23].

\indent By taking into account these point our  aim in this paper
is study of the spectrum and also to solve the inverse problem for
singular non-self-adjoint operators by transmission coefficient
and normalizing numbers corresponding to quasi-eigenfunctions of
the Sturm- Liouville operator with complex periodic potential and
 discontinuous coefficients on the axis.

 1 We will investigate  the representation of fundamental
solutions and  some  relations between  these  solutions which
will be very useful for solving of an inverse problem. In our
investigations the basic role will be played  a good
representation of the solutions of the basic equation (1). As the
potential allows bounded analytic continuation to the upper half
of the complex plane , we can conduct detailed analysis of problem
(2)-(5).

 2 We will study the spectrum of the problem (2)-(5).
  As our solutions have an obvious form that kernel of the resolvent
   will have a simple representation and consequently it will be possible
   to study its analytical properties in detail.

3  We will put and solve the inverse problem.  As it was pointed
out in the introduction some of the authors for solving of the
inverse problems employed operator transformations and studied
solvability of the obtained Gelfand-Levitan type integral
equations. In our case, to solve the inverse problem we have to
solve a system of algebraic equations. So inverse problem has a
unique solution and the numbers $q_{n} $,$p_{n} $  and  $\beta $
are defined constructively by the spectral data.

  \begin{center}
\noindent \textbf{\underbar{ REPRESENTATION OF FUNDAMENTAL
 SOLUTIONS.}}
\end{center}
  Here we study the solutions of the main equation

\noindent

\[-y''(x)+2\lambda p(x)y(x)+q(x)y(x)=\lambda ^{2} \rho \left(x\right)y(x)\]
that will be convenient in future. We can prove the existence of
these solutions if the condition
\eqref{GrindEQ__3_}-\eqref{GrindEQ__5_} is fulfilled for the
potential. This will be unique restriction on the potential and
later on we'll consider it to be fulfilled.

\noindent \textbf{Theorem 1.  } Let \textit{$q(x)$}and $p(x)$ be
of the form \eqref{GrindEQ__3_}, \eqref{GrindEQ__4_} and
\textit{$\rho \left(x\right)$ }satisfy condition
\eqref{GrindEQ__5_}. Then equation $Ly=\lambda ^{2} y$ has special
solutions of the form

\begin{equation} \label{GrindEQ__6_}
f_{1}^{\pm} (x,\lambda )=e^{ \pm i \lambda x} \left(1+\sum
_{n=1}^{\infty }v_{n}^{\pm} e^{inx} +\sum _{n=1}^{\infty }\sum
_{\alpha =n}^{\infty }\frac{v_{n\alpha }^{\pm} }{n\pm2\lambda }
e^{i\alpha x}    \right)\,\,\, for \,\, x \geq 0
\end{equation}

and

\begin{equation} \label{GrindEQ__7_}
f_{2}^{\pm} (x,\lambda )=e^{\pm\lambda \beta x} \left(1+\sum
_{n=1}^{\infty }v_{n}^{\pm} e^{inx} +\sum _{n=1}^{\infty }\sum
_{\alpha =n}^{\infty }\frac{v_{n\alpha }^{\pm} }{n\mp2i\lambda
\beta } e^{i\alpha x}    \right)\,\,\, for\,\,  x<0
\end{equation}

where the numbers \textit{$v_{n}^{\pm} $} and\textit{ $v_{n\alpha
}^{\pm} $} are determined from the following recurrent relations

\begin{equation} \label{GrindEQ__8_}
\alpha ^{2} \nu _{\alpha }^{\pm} +\alpha \sum _{n=1}^{\alpha }\nu
_{n\alpha }^{\pm} +\sum _{s=1}^{\alpha -1}\left(q_{\alpha -s} \nu
_{s}^{\pm} -p_{\alpha -s} \sum _{n=1}^{s}\nu _{ns}^{\pm}  \right)
+q_{\alpha } =0 ,
\end{equation}

\begin{equation} \label{GrindEQ__9_}
\alpha (\alpha -n)\nu _{n\alpha }^{\pm} +\sum _{s=n}^{\alpha
-1}(q_{\alpha -s} \mp n\cdot p_{\alpha -s} )\nu _{ns}^{\pm} =0 ,
\end{equation}

\begin{equation} \label{GrindEQ__10_}
\alpha \nu _{\alpha }^{\pm} \pm \sum _{s=1}^{\alpha -1}\nu
_{s}^{\pm} p_{\alpha -s} \pm p_{\alpha } =0 .
\end{equation}
and the series (6,7) admits double termwise differentiation.

\noindent    \textbf{Remark1:  } If $\lambda \ne -\frac{n}{2} $
and $Im\lambda > 0$, then $f_{1}^{+} \left(x,\lambda \right)\in
L_{2} \left(0,+\infty \right)$. \textit{}

\noindent \textbf{Remark2:}  If $\lambda \ne -\frac{in}{2\beta } $
and $Re\lambda > 0$, then $f_{2}^{+} \left(x,\lambda \right)\in
L_{2} \left(-\infty ,0\right)$.

By analogy [22], it is easy to see that equation
\eqref{GrindEQ__1_} for  $\lambda\neq0$,$\lambda\neq
\pm\frac{n}{2}$ è $ \lambda\neq \pm\frac{in}{2\beta}$ has
fundamental solutions $f_{1}^{+} \left(x,\lambda
\right)$,$f_{1}^{-} \left(x,\lambda \right)$ ($f_{2}^{+}
\left(x,\lambda \right)$,$f_{2}^{-} \left(x,\lambda \right)$) for
which

\[W\left[f_{1}^{+} (x,\lambda ),f_{1}^{-} (x,\lambda )\right]=2i\lambda ,\]

\[W\left[f_{2}^{+} (x,\lambda ),f_{2}^{-} (x,\lambda )\right]=-2i\lambda \beta \]
 (where $W\left[f,g\right]=f'g-fg'$) is satisfied\textbf{}

Then use the conjunction conditions

\[ y\left(0+\right)=y(0-),\]
\[ y'\left(0+\right)=y'(0-)\]

we get  that each solution of equation \eqref{GrindEQ__1_} may be
represented as linear combinations of these solutions
\begin{equation}\label{11}
    f_{2}^{+}\left(x,\lambda \right)=C_{11} \left(\lambda \right)f_{1}^{+} \left(x,\lambda \right)+C_{12} \left(\lambda \right)f_{1}^{-} \left(x,\lambda \right), \, \, \, for\, \, \, x\ge
    0,
\end{equation}

\begin{equation}\label{12}
    f_{1}^{+} \left(x,\lambda \right)=C_{22} \left(\lambda \right)f_{2}^{+} \left(x,\lambda \right)+C_{21} \left(\lambda \right)f_{2}^{-} \left(x,\lambda \right)\, ,\, \, \, \, \, \, \, for\, \, \, x<0\, \, ,
\end{equation}

were

\[C_{11} \left(\lambda \right)=\frac{W[f_{2}^{+} \left(0,\lambda
\right),f_{1}^{-}\left(0,\lambda \right)]}{2i\lambda },\]

\[C_{12} \left(\lambda \right)=\frac{W[f_{1}^{+} \left(0,\lambda \right),f_{2}^{+} \left(0,\lambda \right)]}{2i\lambda } ,\]

\begin{equation} \label{GrindEQ__13_}
C_{22} \left(\lambda \right)=\frac{W[f_{1}^{+} \left(0,\lambda
\right),f_{2}^{-} \left(0,\lambda \right)]}{2\beta\lambda } ,
\end{equation}
\begin{equation} \label{GrindEQ__14_}
C_{21} \left(\lambda \right)=\frac{W[f_{2}^{+} \left(0,\lambda
\right),f_{1}^{+} \left(0,\lambda \right)]}{2\beta\lambda } .
\end{equation}

 By  (13) è (14) the coefficients  $C_{22}(\lambda)$ ,
 $C_{21}(\lambda)$  can be  predetermined and at  $ x\geq0$  we get

\begin{equation} \label{GrindEQ__15_}
C_{22} \left(\lambda \right)=\frac{i}{\beta } C_{11}
\left(-\lambda \right),  C_{21} \left(\lambda
\right)=-\frac{i}{\beta } C_{12} \left(\lambda \right).
\end{equation}

\indent Let
\begin{equation} \label{GrindEQ__16_}
f_{n}^{\pm } (x)=\mathop{\lim }\limits_{\lambda \to \mp
\frac{n}{2} } (n\pm 2\lambda )f_{1}^{\pm}(x,\lambda )=\sum
_{\alpha =n}^{\infty }v_{n\alpha }^{\pm} e^{i\alpha x}
e^{-i\frac{n}{2}x} ,
\end{equation}

It follows from relation \eqref{GrindEQ__6_} that  $f_{n}^{\pm }
(x)\not  \equiv 0$ is valid for $v_{n\alpha }^{\pm} \ne 0$.

From this we obtained that
 $W[f_{n}^{\pm } \left(x\right),f_{1}^{\mp }
\left(x,\mp \frac{n}{2} \right)]=0$ and consequently the functions
$f_{n}^{\pm } \left(x\right),f_{1}^{\mp } \left(x,\mp \frac{n}{2}
\right)$ are linear dependent.

 Therefore

 \begin{equation} \label{GrindEQ__17_}
f_{n}^{\pm } \left(x\right)=S_{n}^{\pm} f_{1}^{\mp} \left(x,
\frac{n}{2} \right).
\end{equation}
 Comparing these formulas we get that  $S_{n}^{\pm}=v_{nn }^{\pm}.$

\begin{center}
\noindent \textbf{\underbar{ SPECTRUM OF THE OPERATOR L }}
\end{center}

\noindent To study the spectrums of the operator $L$ at first we
calculate the kernel of the resolvent of the operator
\textbf{$\left(L-\lambda ^{2} I\right)$ }. Divided  the plane
$\lambda$ into sectors
\[S_{k} =\{ {\raise0.7ex\hbox{$ k\pi  $}\!\mathord{\left/{\vphantom{k\pi  2}}\right.\kern-\nulldelimiterspace}\!\lower0.7ex\hbox{$ 2 $}} <\arg \lambda <{\raise0.7ex\hbox{$ (k+1)\pi  $}\!\mathord{\left/{\vphantom{(k+1)\pi  2}}\right.\kern-\nulldelimiterspace}\!\lower0.7ex\hbox{$ 2 $}} \} ,k=\overline{0,3}.\]

 By means of general methods for the kernel of the resolvent, we
 get

 \begin{equation} \label{GrindEQ__18_}
R_{11} \left(x,t,\lambda \right)=\frac{1}{W[f_{1}^{+} ,f_{2}^{+}
]} \left\{\begin{array}{c} {f_{1}^{+} \left(x,\lambda
\right)f_{2}^{+} \left(t,\lambda \right)\, \, \, \, \, \, \, \, \,
for\, \, t<x} \\ {f_{1}^{+} \left(t,\lambda \right)f_{2}^{+}
\left(x,\lambda \right)\, \, \, \, \, \, \, \, \, for\, \, t>x}
\end{array}\right.        \lambda \in S_{0} ,
\end{equation}

\begin{equation} \label{GrindEQ__19_}
R_{12} \left(x,t,\lambda \right)=\frac{1}{W[f_{1}^{+} ,f_{2}^{-}
]} \left\{\begin{array}{c} {f_{1}^{+} \left(x,\lambda
\right)f_{2}^{-} \left(t,\lambda \right)\, \, \, \, \, \, \, \, \,
for\, \, t<x} \\ {f_{1}^{+} \left(t,\lambda \right)f_{2}^{-}
\left(x,\lambda \right)\, \, \, \, \, \, \, \, \, for\, \, t>x}
\end{array}\right.        \lambda \in S_{1} ,
\end{equation}

\begin{equation} \label{GrindEQ__20_}
R_{21} \left(x,t,\lambda \right)=\frac{1}{W[f_{1}^{-} ,f_{2}^{-}
]} \left\{\begin{array}{c} {f_{1}^{-} \left(x,\lambda
\right)f_{2}^{-} \left(t,\lambda \right)\, \, \, \, \, \, \, \, \,
for\, \, t<x} \\ {f_{1}^{-} \left(t,\lambda \right)f_{2}^{-}
\left(x,\lambda \right)\, \, \, \, \, \, \, \, \, for\, \, t>x}
\end{array}\right.        \lambda \in S_{2} ,
\end{equation}

\begin{equation} \label{GrindEQ__21_}
R_{22} \left(x,t,\lambda \right)=\frac{1}{W[f_{1}^{-} ,f_{2}^{+}
]} \left\{\begin{array}{c} {f_{1}^{-} \left(x,\lambda
\right)f_{2}^{+} \left(t,\lambda \right)\, \, \, \, \, \, \, \, \,
for\, t<x} \\ {f_{1}^{-} \left(t,\lambda \right)f_{2}^{+}
\left(x,\lambda \right)\, \, \, \, \, \, \, \, for\, t>x}
\end{array}\right.        \lambda \in S_{3} .
\end{equation}

\indent \textbf{\underbar{Remark 3 :}} $ R_{\lambda  }  = (L -
\lambda ^2 I)^{-1}$ exists and bounded for all $ \lambda $ out of
$\{Re\lambda=0\}\bigcup\{Im\lambda=0\}$  and $ C_{11} \left( { \pm
\lambda } \right) \ne 0$, $ C_{12} \left( { \pm \lambda } \right)
\ne 0$.

\textbf{\underbar{Theorem 2.}   }$L$ has no eigenvalues for real
and pure imaginary $\lambda $. It's continuous  spectra consist of
axes $Re\lambda =0$  and $Im\lambda =0$ on which there may exist
spectral singularities coinciding with the numbers
$\frac{in}{2\beta } ,\, \, \, \frac{n}{2} ,\, \, \, n=\pm 1,\pm
2,\pm 3,...$

 \textbf{\underbar{Proof:}}

     We note that equation \eqref{GrindEQ__1_}has fundamental solutions$f_{1}^{+} \left(x,\lambda
\right)$,$f_{1}^{-} \left(x,\lambda \right)$ ($f_{2}^{+}
\left(x,\lambda \right)$,$f_{2}^{-} \left(x,\lambda \right)$)
on$\left|Im\lambda \right|<\frac{\varepsilon }{2} \, \, \, \,
(\left|Re\lambda \right|<\frac{\varepsilon }{2} )$.

\noindent Then for $Im\lambda =0$  solution of equation
\eqref{GrindEQ__1_} can be written in the form

\[y\left(x,\lambda \right)=C_{1} f_{1}^{+} \left(x,\lambda \right)+C_{2} f_{1}^{-} \left(x,\lambda \right)\, \, \, \]
In case $Im\lambda =0$ the solution  $f_{1}^{\pm } \left(x,\lambda
\right)$  has the form

\[f_{1}^{\pm } (x,\lambda )=e^{\pm iRe\lambda x} \left(1+\sum
_{n=1}^{\infty }v_{n}^{\pm} e^{inx}+ \sum _{n=1}^{\infty
}\frac{1}{n\pm 2\lambda } \sum _{\alpha =n}^{\infty }V_{n\alpha }
e^{i\alpha x}   \right)\, \, \, \, \, \, \, \, \] then
$y\left(x,\lambda \right)\in L_{2} \left(-\infty ,+\infty \right)$
except when $C_{1} =C_{2} =0$.

\noindent Analogously we can prove case $Re\lambda =0$.

 Since in $\left|Re\lambda \right|<\frac{\varepsilon }{2} $  the
functions$f_{2}^{+} \left(x,\lambda \right)$, $f_{2}^{-}
\left(x,\lambda \right)$ form fundamentals solutions, then

\[y\left(x,\lambda \right)=C_{3} f_{2}^{+} \left(x,\lambda \right)+C_{4} f_{2}^{-} \left(x,\lambda \right)\, \, \, .\]
If  $Re\lambda =0$ then the solution $f_{2}^{\pm } \left(x,\lambda
\right)$ has the form

\[f_{2}^{\pm } (x,\lambda )=e^{\pm iIm\lambda \beta x} \left(1+\sum
_{n=1}^{\infty }v_{n}^{\pm} e^{inx}+\sum _{n=1}^{\infty
}\frac{1}{n\pm 2Im\lambda \beta } \sum _{\alpha =n}^{\infty
}V_{n\alpha } e^{i\alpha x}   \right)\, \] then  $y\left(x,\lambda
\right)\in L_{2} \left(-\infty ,+\infty \right)$  except when
$C_{3} =C_{4} =0$.

 Taking into account Remark 2, we'll investigate the
function

\[
R(x,t,\lambda ) = \left\{ {\begin{array}{*{20}c}
   {R_{11} (x,t,\lambda )}  \\
   {R_{12} (x,t,\lambda )}  \\
   {R_{21} (x,t,\lambda )}  \\
   {R_{22} (x,t,\lambda )}  \\
\end{array}} \right.\begin{array}{*{20}c}
   ,  \\
   ,  \\
   ,  \\
   ,  \\
\end{array}\begin{array}{*{20}c}
   {\lambda  \in S_0 }  \\
   {\lambda  \in S_1 }  \\
   {\lambda  \in S_2 } \\
   {\lambda  \in S_3 } \\
\end{array}
\]

 in  the neighborhood of poles $\lambda^{\ast}$ from  $\{Re \lambda = 0\} \bigcup \{Im \lambda = 0\}$.
  Without loss of generality we can investigate $R(x,t,\lambda)$
  in the neighborhood of poles $\lambda^{\ast}_{0}$ from $[0,\infty)$.
   Then  the number $\lambda^{\ast}_{0}$ coincides with one of the numbers
    $ \left( {\frac{n}{{2}}} \right)$ , $n\in N $.  From
\eqref{GrindEQ__18_}-\eqref{GrindEQ__21_} it follows that the
limit $ \mathop {\lim }\limits_{\lambda  \to \lambda^{\ast}_{0} }
(\lambda - \lambda^{\ast}_{0} )R\left( {x,t,\lambda } \right) =
R_0 \left( {x,t} \right) $ exists and $ R_0 \left( {x,t} \right)$
is a bounded function with respect to all the variables. Let
$\theta(x)$ be an arbitrary finite function. Then $ \varphi \left(
x \right) = \int\limits_{ - \infty }^{ + \infty } {R_0 \left(
{x,t} \right)} \theta \left( t \right)dt$ is a bounded solution of
equation (1) for $ \lambda=\lambda^{\ast}_{0}$. Therefore  $
\varphi \left( x \right) = C_0 f_1^ +  \left(
{x,\lambda^{\ast}_{0} }\right)$ . Comparison of the last relation
with formulae \eqref{GrindEQ__18_}-\eqref{GrindEQ__21_} shows that
if $ \lambda^{\ast}_{0} \ne \frac{n}{2} ,n \in N$ then $ C_0=0$
and so the kernel of the resolvent has removable singularity at
the point $ \lambda^{\ast}_{0}$. So there is a case
$\lambda^{\ast}_{0} = \frac{n}{2} ,n \in N$ where the kernel of
resolvent has poles of the first order. Since $ f_1^ + \left(
{x,\lambda^{\ast}_{0} } \right) \notin L_2 \left( { - \infty , +
\infty } \right)$  then $ \lambda^{\ast}_{0}=(\frac{n}{2})$ is a
spectral singularity of the operator $L$ in the sense of
M.N.Naimark [24, p.450].( Analogously we can show that the
resolvent has poles of the first order at the points $
\frac{in}{2\beta } ,\, \, \, \frac{n}{2} ,\, \, \, n=\pm 1,\pm
2,\pm 3,...$ which are spectral singularities of the operator
$L$).

 \noindent In order to show that all numbers from $\{Re \lambda = 0\} \bigcup \{Im \lambda = 0\}$
 belong to the continuous spectra we introduce the Dirichlet  condition at $x=0$, the resulting operator
 $L^{0}$ splits into the direct sum of corresponding Dirichlet operators $L^{-}$, $L^{+}$ on the negative and positive
 semi-axes. Therefore the essential spectrum of $L^{0}$ is
  $\{Re \lambda = 0\} \bigcup \{Im \lambda = 0\}$, and the initial operator $L$ has the same essential
 spectrum. Since eigenvalues of infinite multiplicities are not
 possible, the essential spectrum of $L$ coincide with the
 continuous one.
  The Theorem  is proved.

\noindent \textbf{Lemma 1.}   The eigenvalues of operator $L$ are
finite and coincide with zeros of the functions $C_{12} (\lambda
),\, \, C_{11} (-\lambda ),\, \, C_{12} (-\lambda ),\, \, C_{11}
(\lambda )$ from sectors  $S_{k} =\{ {\raise0.7ex\hbox{$ k\pi
$}\!\mathord{\left/{\vphantom{k\pi
2}}\right.\kern-\nulldelimiterspace}\!\lower0.7ex\hbox{$ 2 $}}
<\arg \lambda <{\raise0.7ex\hbox{$ (k+1)\pi
$}\!\mathord{\left/{\vphantom{(k+1)\pi
2}}\right.\kern-\nulldelimiterspace}\!\lower0.7ex\hbox{$ 2 $}} \}
,k=\overline{0,3}$ respectively.

\noindent \textbf{Definition 1.}  The data $\{ \,  C_{11}
\left(\lambda \right),\, \, C_{12} \left(\lambda \right)\} $ are
called the spectral data of $L$.

\begin{center}
\noindent \textbf{\underbar{ INVERSE PROBLEM}}
\end{center}

 From the representation \eqref{GrindEQ__18_}-\eqref{GrindEQ__21_} it also
follows that for each $x$ and $t$ from $ (-\infty,+\infty)$ the
kernel $ R(x,t,\lambda)$ admits a meromorphic continuation from
the sectors $ S_{k}, k=0,1,2,3$ and may have poles at the points $
\left( {\pm\frac{n}{{2}}} \right)$ and $\left(
{\pm\frac{in}{{2\beta}}} \right)$, $n\in N $, outside of $S_{k},
k=0,1,2,3$. Such poles of the resolvent are called
quasi-stationary states of the operator $L$. Thus the
quasi-stationary states of the operator $L$ are the numbers $
\left( {\pm\frac{n}{{2}}} \right)$ and $\left(
{\pm\frac{in}{{2\beta}}} \right)$, $n\in N $. The numbers
$V_{nn}$, $n \in N$ play a part of normalizing numbers
corresponding to quasi-eigenfunction of the operator $L$. So, it
is make natural the formulation of the inverse problem about
reconstruction of the potential of the equation
\eqref{GrindEQ__1_} and the number $\beta$.

Let's study the inverse problem for the problem (1-3). The inverse
problem is formulated as follows.

 Given the spectral data
  $\{\, \, C_{11} \left(\lambda \right),\, \, C_{12} \left(\lambda
\right)\} $ construct the potential
$q\left(x\right)$,$p\left(x\right)$ and $\beta $.

Using the results obtained above we arrive at the following
procedure for solution of the inverse problem.

\noindent 1. Taking into account (16) it is easy to check that

\[\mathop{\lim }\limits_{\lambda \to {\raise0.7ex\hbox{$ n $}\!\mathord{\left/{\vphantom{n 2}}\right.\kern-\nulldelimiterspace}\!\lower0.7ex\hbox{$ 2 $}} } \left(n-2\lambda \right)\frac{C_{11} \left(\lambda \right)}{C_{12} \left(\lambda \right)} =v_{nn}^{-} ,\mathop{\lim }\limits_{\lambda \to -{\raise0.7ex\hbox{$ n $}\!\mathord{\left/{\vphantom{n 2}}\right.\kern-\nulldelimiterspace}\!\lower0.7ex\hbox{$ 2 $}} } \left(n+2\lambda \right)\frac{C_{12} \left(\lambda \right)}{C_{11} \left(\lambda \right)} =v_{nn}^{+} \]
consequently we find all numbers $v_{nn}^{\pm } ,\, \, n=1,2,...$.

\noindent 2. Taking into account (16) we get

\[\sum _{\alpha =B}^{\infty }v_{n\alpha }^{\pm } e^{i\alpha x} e^{-i\frac{n}{2} x} =v_{nn}^{\pm } \cdot e^{i\frac{n}{2} x} \left(1+\sum _{n=1}^{\infty }v_{n}^{\mp } e^{inx} +\sum _{n=1}^{\infty }\sum _{\alpha =n}^{\infty }\frac{v_{n\alpha }^{\mp } }{n+n} e^{i\alpha x}    \right) \]
or

\[v_{m,\alpha +m}^{\pm } =v_{mm}^{\pm } \left(v_{\alpha }^{\mp } +\sum _{n==1}^{\alpha }\frac{v_{n\alpha }^{\mp } }{n+m}  \right)\]
from which all numbers  $v_{n\alpha }^{\pm } ,\, \, \, \alpha
=1,2.....,\, \, \, n=1,2,....n<\alpha $ are defined.

\noindent 3. Then from recurrent formula (8-10), find all numbers
$p_{n} $ and $q_{n} $.

\noindent 4. The number $\beta $  is defined from equality

\[\beta =iC_{11} \left(\lambda _{n} \right)C_{11} \left(-\lambda _{n} \right).\]
Really from (11-12) we derive for $\lambda _{n} \in S_{0} $ which
is a root of the equation $C_{12}(\lambda)=0$ the true relation

\noindent $C_{11} \left(\lambda _{n} \right)=\frac{f_{2}^{+}
\left(x,\lambda _{n} \right)}{f_{1}^{+} \left(x,\lambda _{n}
\right)} $, $C_{22} \left(\lambda _{n} \right)=\frac{f_{1}^{+}
\left(x,\lambda _{n} \right)}{f_{2}^{+} \left(x,\lambda _{n}
\right)} $ i.e. $C_{11} \left(\lambda _{n} \right)C_{22}
\left(\lambda _{n} \right)=1$.

\noindent Then  we get $\beta =iC_{11} \left(\lambda _{n}
\right)C_{11} \left(-\lambda _{n} \right)$.

So inverse problem has a unique solution and the numbers $q_{n}
$,$p_{n} $  and  $\beta $  are defined constructively by the
spectral data.

\indent \textbf{Theorem 3.} The specification of the spectral data
uniquely determines potentials $p\left(x\right)$,
$q\left(x\right)$ and $\beta $.

\noindent \textbf{REFERENCES.}

\noindent 1 T. Aktosun, M. Klaus, and C. van der Mee, Wave
scattering in one dimension with absorption, J. Math. Phys. 39,
1957-1992, 1998.

 \noindent 2 T. Aktosun, M. Klaus, and C. van der Mee, Inverse
scattering in one-dimensional nonconservative media, Integral
Equations and Operator Theory 30, 297-316 ,1998

 \noindent 3. Jaulent M. and C. Jean, The inverse s-wave scattering problem
for a class of potentials depending on energy, Commun. Math. Phys.
28, 177-220 ,1972.

 \noindent 4. Jaulent M. and C. Jean, The inverse problem for the
one-dimensional Schrodinger equation with an energy-dependent
potential. I, Ann. Inst. Henri Poincare´, Sect. A 25, 105-118
,1976.

  \noindent 5. Jaulent M and C. Jean, The inverse problem for the
one-dimensional Schrodinger equation with an energy-dependent
potential. II,  Ann. Inst. Henri Poincare´, Sect. A 25, 119-137,
1976.

 \noindent 6. Jaulent M., Inverse scattering problems in absorbing media, J.
Math. Phys. 17, 1351-1360 , 1976.

 \noindent 7. Pivovarchik V.N. Reconstruction of the potential of the
Sturm-Liouville equation from three spectra of - Functional
Analysis and Its Applications, 1999

\noindent 8.  Freiling G,  Yurko V Inverse problems for
differential equations with turning points   - Inverse problems,
1997

\noindent 9. Eberhard W. and Freiling G.,  The distribution of the
eigenvalues for second order eigenvalue problems in the presence
of an arbitrary number of turning points, Results in Mathematics,
21, 24-41, 1992

\noindent 10. Freiling G. and Yurko V.A. , Inverse Sturm-Liouville
Problems and their Applications, NOVA Science Publishers, New
York, 2001.

\noindent 11.Yurko V   Inverse spectral problems for differential
pencils on the half-line with turning points  - Journal of
Mathematical Analysis and Applications, 2006

\noindent 12 Krein MG On determination of the potential of a
particle from its S-function  - Dokl. Akad. Nauk SSSR, 1955

\noindent 13  Blagoveshchensky A.S.,  A one-dimensional inverse
boundary value problem for a second order hyperbolic equation.
(Russian) Zap. Naucn. Sem. Leningrad. Otdel. Mat. Inst. Steklov.
(LOMI) 15 1969 85--90.

\noindent 14.  Grinberg N. I..1991 The one-dimensional inverse
scattering problem for the wave equation, Math. USSR Sb. 70,
557-572, [Mat. Sb. 181 (8), 1114-1129 (1990) (Russian)].

\noindent 15. Belishev MI  Inverse spectral indefinite problem for
the equation $ y''+ \lambda  p(x) y= 0$ on an interval -
Functional Analysis and Its Applications, 1987

 \noindent 16.  Bethe H.A  A continuum theory of the compound nucleus
 - Physical Review,57,1125, 1940

 \noindent 17.  Feshbach H ,  Porter CE,  Weisskopf VF Model for nuclear
 reactions with neutrons - Physical Review, 1954

\noindent 18.   Gasymov M.G. Spectral analysis of a class
non-self-adjoint operator of the  second order. Functional
analysis and its applications. (In Russian), V34.1.pp.14-19,1980

 \noindent 19. Gasymov MG.,Gusejnov GS  Determination of diffusion operators according to
spectral data- Dokl., Akad. Nauk Az. SSR,(In Russian) 1981

\noindent 20.  Efendiev, R.F. The Characterization Problem for One
Class of Second Order  Operator Pencil with Complex Periodic
Coefficients. Moscow Mathematical   Journal, 2007, Volume 7,Number
1, pp.55-65.

\noindent 21.  Efendiev, R. F. Complete solution of an inverse
problem for one class of the  high order ordinary differential
operators with periodic coefficients. Zh. Mat. Fiz. Anal. Geom.
2006, V.2, no. 1, 73--86,

\noindent 22.  Efendiev, R.F. Spectral analysis of a class of
non-self-adjoint  differential operator pencils with a generalized
function.Teoreticheskaya i  Matematicheskaya  Fizika,
2005,Vol.145,1.pp.102-107,October (in Russian).
 Theoretical and Mathematical Physics,1451(1):1457-461,(English).

\noindent 23. Efendiev, R. F Inverse indefinite spectral problem
for second order  differential operator  with    complex periodic
coefficients- Proceeding of Baku State University,(In Russian),
2008, ¹4, 17-22

 \noindent   24.  Naimark,  M.A. Linear differential operators.Moscow, Nauka, 1969(Russian).
 \end{document}